\documentclass{amsart}

\usepackage[english]{babel}
\usepackage{leftidx}
\usepackage[numbers]{natbib}
\usepackage{bigints}
\usepackage{empheq}
\usepackage{cancel}
\usepackage{array}
\usepackage{graphicx}
\usepackage{lipsum}
\usepackage{bigints}
\usepackage{amssymb}
\usepackage{cancel}
\usepackage{amsthm}
\newenvironment{myproof}[2] {\paragraph{\textbf{Proof of {#1} {#2} :}}}{\hfill$\square$}

\usepackage[table]{xcolor}
\usepackage{needspace}
\usepackage{geometry}
\newcommand{\divergence}[1]{\textup{div}\,#1}
\usepackage{enumitem}
\newlist{myQuoteEnumerate}{enumerate}{1}
\setlist[myQuoteEnumerate,1]{label=(P1)}

\newlist{myQuoteEnumerate1}{enumerate}{1}
\setlist[myQuoteEnumerate1,1]{label=(P2)}

\newtheorem{theorem}{Theorem}[section]

\newtheorem{question-non}[]{}

\newtheorem{cor}[theorem]{Corollary}
\newtheorem{definition}[theorem]{Definition}
\newtheorem{example}[theorem]{Example}

\title{Triviality results for quasi $k$-Yamabe solitons}

\author{Tokura, W. $^{{1},\ast}$}
\address{$^{1}$ Universidade Estadual de Mato Grosso do Sul, 79150-000, Av. João Pedro Fernandes, 2101 - Centro, Maracaju, MS, Brazil.}
\email{williamisaotokura@hotmail.com $^{1}$}

\author{Batista, E. $^{2}$}
\address{$^{2}$ Universidade Federal de Goi\'as, IME, 131, 74001-970, Goi\^ania, GO, Brazil.}
\email{elismardb@gmail.com.br $^{2}$}

\author{Kai, P.  $^{3}$}
\address{Universidade Federal de Goi\'{a}s, INF, s/n, 74690-900, Goi\^{a}nia,
GO, Brazil.}
\email{priscila.kai@hotmail.com $^{3}$}

\thanks{$^{\ast}$ 
Corresponding author}



\begin{document}

\begin{abstract}
In this paper, we show that any compact quasi $k$-Yamabe gradient solitons must have constant $\sigma_{k}$-curvature. Moreover, we provide a certain condition for a compact quasi $k$-Yamabe soliton to be gradient, and for noncompact solitons, we present a geometric rigidity under a decaiment condition on the norm of the soliton field.
\end{abstract}

\maketitle

\section{Introduction and main results}
\label{intro}


A Riemannian manifold $(M^n,g)$, $n\geqslant3$, is an \textit{Einstein-type manifold} if there exists a vector field $X\in \mathfrak{X}(M)$ and a smooth function $\lambda:M\rightarrow \mathbb{R}$ such that
\begin{equation}\label{et}\alpha Ric+\frac{\beta}{2}\mathcal{L}_{X}g+\mu X^{\flat}\otimes X^{\flat}=(\rho R+\lambda)g,
\end{equation}
for some constants $\alpha$, $\beta$, $\mu$ $\rho\in\mathbb{R}$ with $(\alpha,\beta,\mu)\neq(0,0,0)$. Here $\mathcal{L}_{X}g$, $X^{\flat}$ and $R$ stand, respectively, for the Lie derivative of $g$ in the direction of $X$, the 1-form metrically dual to the vector field $X$ and the scalar curvature of $M$. The concept of Einstein-type manifold was introduced recently by Catino \textit{et al.} \cite{catino2016geometry}. It is worth noting that in terms of equation \eqref{et}, an Einstein-type manifold unifies various particular cases well studied in the literature, such as gradient Ricci solitons, gradient almost Ricci solitons, Yamabe solitons, quasi Yamabe solitons, conformal gradient solitons, quasi Einstein manifolds and $\rho$-Einstein solitons. Each of them has a particular importance.

Our purpose is to study some cases of the Einstein-type manifolds which were not addressed in \cite{catino2016geometry}. More precisely, we focus our analysis on the class $\alpha=0$, $\beta=1$, $\mu=-\frac{1}{m}$ and $\lambda=\rho R-\sigma_{k}+c$, $c\in\mathbb{R}$ where $\sigma_k$ is the $\sigma_k$-curvature of $g$. We recall that, if we denote by $\lambda_{1},\lambda_{2},\dots,\lambda_{n}$ the
eigenvalues of the symmetric endomorphism $g^{-1}A_{g}$, where $A_{g}$ is the Schouten tensor defined by
\begin{equation*}
A_{g}=\frac{1}{n-2}\left[Ric_{g}-\frac{R}{2(n-1)}g\right],
\end{equation*}
then the $\sigma_{k}$-curvature of $g$ is defined as the $k$-th symmetric elementary function of $\lambda_{1},\dots,\lambda_{n}$, namely
\begin{equation*}
\sigma_{k}=\sigma_{k}(g^{-1}A_{g})=\sum_{i_{1}<\dots <i_{k}}\lambda_{i_{1}}\cdot \dots\cdot \lambda_{i_{k}}, \quad \text{for}\quad 1\leq k\leq n.
\end{equation*}
 
In this sense, we introduce the following class of manifolds.

\begin{definition}
A Riemannian manifold $(M^n,g)$, $n\geqslant3$, is a \textit{quasi $k$-Yamabe soliton} if there exists a vector field $X\in \mathfrak{X}(M)$ and two constants $m$, $\lambda$ (where $m$ is not zero) such that
\begin{equation}\label{fundamental2*}\frac{1}{2}\mathcal{L}_{X}g-\frac{1}{m}X^{\flat}\otimes X^{\flat}=(\sigma_{k}-\lambda)g,
\end{equation}
We will write the soliton in \eqref{fundamental2*} as $(M^{n}, g, X, \lambda)$ for the
sake of simplicity. If $X=\nabla f$ for some smooth function $f:M\rightarrow\mathbb{R}$, we say that $(M^n,g,f,\lambda)$ is a \textit{quasi $k$-Yamabe gradient soliton}.  In this case, equation \eqref{fundamental2*} can be rewritten as 
\begin{equation}\label{fundamental*}\nabla^{2}f-\frac{1}{m}df\otimes df=(\sigma_{k}-\lambda)g,
\end{equation}
Moreover, when $f$ is a constant or $X=0$ the soliton is called \textit{trivial}.
\end{definition}

It is worth noting that quasi $k$-Yamabe solitons correspond to a large class of manifolds well studied in the literature. For instance, if $m\to\infty$, then quasi $k$-Yamabe solitons reduces to $k$-Yamabe solitons \cite{barboza2020rigidity, bo2018k, catino2012global,tokura2021triviality}. Also, since $\sigma_{1}=\frac{R}{2(n-1)}$, $1$-Yamabe solitons naturally correspond to Yamabe solitons \cite{brozos2016local,calvino2012three,cao2012structure,chu2013scalar,daskalopoulos2013classification,hamilton1988ricci, hamilton1989lectures, hsu2012note,ma2011properties, tokura2018warped}, and quasi $1$-Yamabe solitons correspond to quasi Yamabe solitons studied in \cite{huang2014classification, wang2013noncompact}. When $X=\nabla f$ is a gradient vector field, quasi $1$-Yamabe gradient solitons correspond to an $f$-almost Yamabe soliton \cite{zeng1060h}.


In recent years, many efforts have been devoted to study the geometry of Yamabe solitons and their generalizations. For instance, Hsu in \cite{hsu2012note} shown that any compact gradient $1$-Yamabe soliton is trivial. For $k>1$, the extension of the previous result was recently investigated. Catino \textit{et al.} \cite{catino2012global} proved that any compact
gradient $k$-Yamabe soliton with nonnegative Ricci curvature is trivial. Bo
\textit{et al.} \cite{bo2018k} also proved that any compact gradient $k$-Yamabe
soliton with negative constant scalar curvature necessarily has constant
$\sigma_{k}$-curvature. In
\cite{tokura2021triviality} it was shown that any compact gradient $k$-Yamabe
soliton must be trivial. For $m\neq\infty$, Huang and Li \cite{huang2014classification} proved that any compact quasi $1$-Yamabe gradient
soliton is trivial. Our first theorem  extends all these results at once.

\begin{theorem}\label{T1}
Any compact quasi $k$-Yamabe gradient soliton $(M^n,g,f,\lambda)$ is trivial, i.e., has constant $\sigma_{k}$-curvature.
\end{theorem}

The Hodge-de Rham decomposition theorem shows that any vector field $X$ on a compact oriented Riemannian manifold $M^n$ can be decompose as follows: 
\begin{equation}\label{Hodge} X = \nabla h + Y,
\end{equation}
where $h:M\rightarrow\mathbb{R}$ is a smooth function and $Y$ is a divergence free vector field on $M^n$. In fact,
by the Hodge-de Rham theorem  \cite{warner2013foundations}, we have that $X^{\flat}$ takes the following form:
\[X^{\flat}=d\alpha+\delta\beta+\gamma.\]
Taking $Y = (\delta\beta+\gamma)^{\sharp}$ and $(d\alpha)^{\sharp}=\nabla h$ we arrive at the desired result.

Our next result provide a necessary and sufficient condition for a compact quasi $k$-Yamabe soliton to be gradient.

\begin{theorem}\label{T2}The compact quasi $k$-Yamabe soliton $(M^n,g,X,\lambda)$ is gradient if, and only if
\begin{equation}\label{tobegrad}
        \begin{split}&\int_{M}Ric(\nabla h,Y)+\frac{1}{m}g(\nabla^{2}h,X^{\flat}\otimes X^{\flat})-\frac{2}{m}g(\nabla^{2}h,dh\otimes dh)+\frac{1}{m^{2}}|dh\otimes dh|^{2} +\frac{n}{2m}(\sigma_{k}-\lambda)|X|^{2}\\
         &\qquad \qquad\qquad \qquad\qquad\qquad\qquad+\frac{2}{m}(\sigma_{k}-\lambda)|\nabla h|^{2}+\frac{3}{2m^2}|X|^{4}\leqslant0.
        \end{split}
        \end{equation}
where $h$ and $Y$ are the Hodge–de Rham decomposition components of $X$.
\end{theorem}

Note that when $m\to\infty$, then Theorem \ref{T2} correspond to Theorem 1.3 of \cite{tokura2021triviality}. On the other hand, combining Theorem \ref{T1} and Theorem \ref{T2}, we have the following result.

\begin{cor}
Any compact quasi $k$-Yamabe soliton $(M^n,g,X,\lambda)$ satisfying \eqref{tobegrad} is trivial, i.e., has constant $\sigma_{k}$-curvature.
\end{cor}



Now we notice that the same result obtained in \cite{aquino2011some} for compact Ricci solitons also works for compact quasi $k$-Yamabe solitons. More precisely, we have the following theorem.

\begin{theorem}\label{T8}Let $(M^n,g,X,\lambda)$ be a compact quasi $k$-Yamabe soliton and   $X=\nabla h+Y$ the Hodge–de Rham decomposition of $X$. If 
\[\int_{M^n}g(\nabla h,X)dv_{g}\leqslant0,\]
then $(M^n,g)$ is trivial, i.e., has constant $\sigma_k$-curvature.

\end{theorem}

The next theorem provide a rigidity result in the scope of noncompact quasi $k$-Yamabe gradient solitons. 

In  \cite{ma2012remarks}, Ma and Miguel study Liouville type theorem of harmonic functions with finite weighted Dirichlet integral and use it to prove a rigidity result for gradient Yamabe solitons with non negative Ricci curvature. We also provide a rigidity result assuming that the weighted integral of the soliton vector field  is finite.

\begin{theorem}\label{T4}
    Let $(M^n,g,f,\lambda)$ be a complete and noncompact quasi $k$-Yamabe gradient 
    soliton satisfying
    \begin{equation*}
        \int_{M^n\setminus B(r)}d(x,x_{0})^{-1}|\nabla f|d\mu<\infty,
    \end{equation*}
    where $d$ is the distance function with respect to $g$, $B(r)$ is
    the open ball of radius $r>0$ centered at $x_{0}$ and $d\mu=e^{-\frac{f}{m}}dv_{g}$ . If $f$ is $f$-subharmonic, then $(M^n,g)$ has constant $\sigma_k$-curvature.
\end{theorem}

Before we present the proof of main results, let's take a look at some examples. For others examples see Section \ref{sectioninvariant}.

\begin{example}\label{exemplo22}
      Identities
    \[
        Ric_{g_{\mathbb{S}^{n}}}=(n-1)g_{\mathbb{S}^{n}},\quad
        R_{g_{\mathbb{S}^{n}}}=n(n-1)\quad\mbox{and}\quad
        A_{g_{\mathbb{S}^n}}=\frac{1}{2}g_{\mathbb{S}^n},
    \]
    rule the Ricci tensor, scalar curvature and Schouten tensor, respectively,
    of the Euclidean sphere $ (\mathbb{S}^n,g_{\mathbb{S}^n}) $. Therefore, we
    have that
    \[
        \sigma_{k}(g_{\mathbb{S}^{n}}^{-1}A_{g_{\mathbb{S}^{n}}})=
        \frac{1}{2^k}\binom{n}{k},\quad1\leqslant k\leqslant n.
    \]
    Then $ (\mathbb{S}^n,g_{\mathbb{S}^n}) $ is a trivial quasi $k$-Yamabe gradient soliton with $\sigma_{k}=\lambda$. 
    
    Note that, according to Theorem \ref{T1}, the Euclidean sphere does not admit any non-trivial structure of quasi $k$-Yamabe gradient soliton.
\end{example}

\begin{example}Consider the hyperbolic space  $\mathbb{H}^{n+1}=\mathbb{R}\times_{e^{t}}\mathbb{R}^{n}$ furnished with the warped product metric \cite{o1983semi}:
\[g=dt^2+e^{2t}g_{\mathbb{R}^n}^{2}.\]
It is
well known that the horospheres of the hyperbolic space are totally umbilical hypersurfaces isometric to $\mathbb{R}^{n}$ and correspond to slices ${\{t_{0}\}}\times \mathbb{R}^{n}$, $t_{0}\in \mathbb{R}$. Hence, taking the inclusion $$i:\{t_{0}\}\times\mathbb{R}^{n}\rightarrow \mathbb{R}\times\mathbb{R}^{n}\quad (t_{0},x)\mapsto ( t_{0},x),$$
we deduce that the height function from the immersion satisfies $f(x)=t_{0}$, and then the standard Euclidean space $\{t_{0}\}\times\mathbb{R}^{n}$ is a trivial quasi $k$-Yamabe gradient soliton immersed into hyperbolic space with potential $f(x)=t_{0}$ and $\lambda=0$.
\end{example}

\begin{example}\label{examplee1}Consider $M^n=\{(x_{1},\dots, x_{n})\in\mathbb{R}^n\;;\; x_{n}>0\}$, $g_{ij}=(k_{0}x_{n})^{-2}\delta_{ij}$, $k_{0}\in(0,\infty)$ and the potential function
\[f(x_{1},\dots,x_{n})=-m\log\left(\frac{k_{1}}{ x_{n}}\right),\qquad k_{1}\in(0,\infty).\]
By a direct computation, we deduce that
\[\sigma_{k}=-\frac{n!}{k!(n-k)!}(-1)^{k-1}\left(\frac{k_{0}^{2}}{2}\right)^{k},\qquad \nabla^{2}f-\frac{1}{m}df\otimes df=-mk_{0}^2g.\]
Then, $(M^n,g)$ is a quasi $k$-Yamabe gradient soliton with 
\[\lambda=-\frac{n!}{k!(n-k)!}(-1)^{k-1}\left(\frac{k_{0}^{2}}{2}\right)^{k}+mk_{0}^2.\]

\end{example}

\begin{example}\label{Example noncompact}Consider the Euclidean subset $M^n=\{(x_{1},\dots, x_{n})\in\mathbb{R}^n\;;\; x_{1}+\dots+x_{n}>0\}$ furnished with the metric tensor $g_{ij}=\delta_{ij}$ and potential function given by

\[f(x_{1},\dots,x_{n})=-m\log\left(x_{1}+\cdots+x_{n}\right).\]
Since $Ric_{g}=0$, $R_{g}=0$, we have that $A_{g}=0$. Consequently $\sigma_{k}=0$ for all $k\in\{1,\cdots, n\}$. On the other hand,  from the potential function $f$, we deduce that $\nabla^{2}f=\frac{1}{m}df\otimes df$. So, $(M^n,g)$ is a quasi $k$-Yamabe gradient soliton with $\lambda=0$.
\end{example}

Example \ref{Example noncompact} shows that for a noncompact quasi $k$-Yamabe gradient soliton, the condition of constant $\sigma_k$-curvature does not imply that the potential function $f$ is constant.

\begin{example}\label{examplee2}Consider $M^n=\{(x_{1},\dots, x_{n})\in\mathbb{R}^n\;;\; x_{1}+\dots+x_{n}>0\}$, $g_{ij}=(x_{1}+\dots+x_{n})^{-2}\delta_{ij}$ and the potential function
\[f(x_{1},\dots,x_{n})=-m\log\left(\frac{1}{ x_{1}+\dots+x_{n}}\right).\]
By a direct computation, we deduce that
\[\sigma_{k}=-\frac{n!}{k!(n-k)!}(-1)^{k-1}\left(\frac{n}{2}\right)^{k},\qquad \nabla^{2}f-\frac{1}{m}df\otimes df=-mng.\]
Then, $(M^n,g)$ is a quasi $k$-Yamabe gradient soliton with 
\[\lambda=-\frac{n!}{k!(n-k)!}(-1)^{k-1}\left(\frac{n}{2}\right)^{k}+mn.\]

\end{example}

Example \ref{examplee1}, Example \ref{Example noncompact} and Example \ref{examplee2} show that the compactness of the manifold can not be discarded in the prove of  Theorem \ref{T1}.





\section{Proofs}
\begin{myproof}{Theorem}{\ref{T1}}
    If $k=1$, then $(M^n,g)$ is a quasi Yamabe gradient soliton and the result is well known from \cite{huang2014classification}. Now, consider $k\geq2$ and suppose by contradiction that $f$ is nonconstant. Set $u=e^{-\frac{f}{m}}$. Then
    \[\nabla u=-\frac{u}{m}\nabla f,\qquad \nabla^{2}u=\frac{u}{m^2}df\otimes df-\frac{u}{m}\nabla^{2}f,\]
    and \eqref{fundamental*} can be rewritten as follows:
    \begin{equation}\label{1111}\nabla^2 u=-\frac{u}{m}(\sigma_{k}-\lambda)g.
    \end{equation}
    Since $f$ is nonconstant, $u$ is also nonconstant, so from Theorem 1.1 of \cite{catino2012global}, we obtain that $(M^n,g)$ is rotationally symmetric and $M^{n}\setminus \{N,S\}$ is locally conformally flat. Here $N,S$  corresponds to the extremal points of  $u$ in $M$. From  \eqref{1111}, we know that $\nabla u$ is a conformal Killing vector field; hence, we can apply Theorem 5.2 of \cite{viaclovsky2000some} (see also Theorem 1 of \cite{han2006kazdan}) to deduce
\begin{equation}\label{1221}
    0=\int_{M^{n}\setminus\{N,S\}}g(\nabla\sigma_{k},\nabla u)dv_{g}=\int_{M^n}g(\nabla\sigma_{k},\nabla u)dv_{g}=\frac{n}{m}\int_{M^n}u\sigma_{k}(\sigma_{k}-\lambda)dv_{g},
\end{equation}
where in the last equality we have used the divergence theorem. On the other hand, again from the divergence theorem, we get
\begin{equation}\label{808080}0=\int_{M^n}\Delta u dv_{g}=-\frac{n}{m}\int_{M^n}u(\sigma_{k}-\lambda)dv_{g}.
\end{equation}
Combining equations \eqref{1221} and \eqref{808080} we arrive at
\[\frac{n}{m}\int_{M^n}u(\sigma_{k}-\lambda)^2 dv_{g}=0,\]
which implies that $\sigma_{k}=\lambda$ and $u$ is harmonic. Since $M^n$ is compact, $u$ is a constant, which leads to a contradiction. This proves that $f$ is constant.

\end{myproof}

\begin{myproof}{Theorem}{\ref{T2}}
From the Hodge-de Rham decomposition $X=\nabla h+Y$, we deduce that
\begin{equation}\label{qgqg}\frac{1}{2}\mathcal{L}_{Y}g=\frac{1}{2}\mathcal{L}_{X}g-\frac{1}{2}\mathcal{L}_{\nabla h}g,
\end{equation}
and
\begin{equation}\label{t11111}
T_{m}:=\frac{1}{2}\mathcal{L}_{Y}g-\frac{1}{m}X^{\flat}\otimes Y^{\flat}-\frac{1}{m}Y^{\flat}\otimes dh=(\sigma_{k}-\lambda)g-\nabla^{2}h+\frac{1}{m}dh\otimes dh.
\end{equation}
Therefore, to prove that $(M^n,g)$ admits a quasi $k$-Yamabe gradient soliton structure, it is necessary and sufficient to show that $T_{m}=0$. From \eqref{qgqg} we arrive at

\begin{equation}\label{norma}
        \begin{split}\int_{M}|T_{m}|^2&=\int_{M}n(\sigma_{k}-\lambda)^{2}-2(\sigma_{k}-\lambda)\Delta h+2(\sigma_{k}-\lambda)\frac{|\nabla h|^{2}}{m}+|\nabla^{2}h|^{2}-\frac{2}{m}g(\nabla^{2}h,dh\otimes dh)\\
        &\qquad +\frac{1}{m^{2}}|dh\otimes dh|^{2}\\
        &=\int_{M}|\nabla^{2}h|^{2}-n(\sigma_{k}-\lambda)^{2}+2(\sigma_{k}-\lambda)\frac{|\nabla h|^{2}-|X|^{2}}{m}-\frac{2}{m}g(\nabla^{2}h,dh\otimes dh)\\
        &\qquad +\frac{1}{m^{2}}|dh\otimes dh|^{2}.
        \end{split}
        \end{equation}
We are going to compute the right-hand side of above equation using the following integral identity 
\begin{equation}\label{tt1}\int_{M^n}2Ric(\nabla h,Y)dv_{g}=\int_{M^n}\left[Ric(X,X)-Ric(\nabla h,\nabla h)-Ric(Y,Y)\right]dv_{g}.
\end{equation}

Taking the divergence of \eqref{qgqg}, we get
\begin{equation}\label{t2}
\begin{split}
\frac{1}{2}div(\mathcal{L}_{Y}g)(Y)&=\frac{1}{2}div(\mathcal{L}_{X}g)(Y)-\frac{1}{2}div(\mathcal{L}_{\nabla h}g)(Y)\\
&=div(\sigma_{k}-\lambda)(Y)+\frac{1}{m}div(X^{\flat}\otimes X^{\flat})(Y)-\frac{1}{2}div(\mathcal{L}_{\nabla h}g)(Y)\\
&=g(\nabla \sigma_{k},Y)+\frac{1}{m}div(X^{\flat}\otimes X^{\flat})(Y)-\frac{1}{2}div(\mathcal{L}_{\nabla h}g)(Y).
\end{split}
\end{equation}
From the Bochner formula (see Lemma 2.1 of \cite{petersen2009rigidity}), we can express \eqref{t2} as follows
\[\frac{1}{2}\Delta |Y|^{2}-|\nabla Y|^{2}+Ric(Y,Y)=2g(\nabla \sigma_{k},Y)+\frac{2}{m}div(X^{\flat}\otimes X^{\flat})(Y)-2Ric(Y,\nabla h)-2g(\nabla \Delta h,Y),\]
and using the compactness of $M^n$, we arrive at
equation
\begin{equation}\label{ddd1}\int_{M}Ric(Y,Y)=\int_{M}\frac{2}{m}div(X^{\flat}\otimes X^{\flat})(Y)-2Ric(Y,\nabla h)+|\nabla Y|^{2}.
\end{equation}
On the other hand, the same argument as above shows that

\begin{equation}\label{113344}\frac{1}{2}\Delta |X|^{2}-|\nabla X|^{2}+Ric(X,X)+(n-2)g(\nabla\sigma_{k},X)+\frac{1}{m}\nabla_{X}
|X|^{2}=\frac{2}{m}div(X^{\flat}\otimes X^{\flat})(X).
\end{equation}
Since
\begin{equation*}
    \begin{split}
        \int_{M^n}|\nabla X|^2dv_{g}&=\int_{M^n}\big{[}|\nabla^2 h|^2+|\nabla Y|^2+2g(\nabla\nabla h,\nabla Y)\big{]}dv_{g}\\
        &=\int_{M^n}\big{[}|\nabla^2 h|^2+|\nabla Y|^2-2g(\nabla\Delta h+Ric(\nabla h), Y)\big{]}dv_{g}\\
        &=\int_{M^n}\big{[}|\nabla^2 h|^2+|\nabla Y|^2-2Ric(\nabla h, Y)\big{]}dv_{g},\\
    \end{split}
\end{equation*}
we may integrate \eqref{113344} over $M^n$ to deduce
\begin{equation}\label{dddd4}
    \begin{split}
        \int_{M^n}Ric(X,X)&=\int_{M^n}|\nabla X|^{2}+\frac{2}{m}div(X^{\flat}\otimes X^{\flat})(X)-(n-2)g(\nabla\sigma_{k},X)-\frac{1}{m}\nabla_{X}
|X|^{2}\\
        &=\int_{M^n}|\nabla^2 h|^2+|\nabla Y|^2-2Ric(\nabla h, Y)+\frac{2}{m}div(X^{\flat}\otimes X^{\flat})(X)\\
        &\qquad  +(n-2)\sigma_{k}\left[(\sigma_{k}-\lambda)n+\frac{1}{m}|X|^{2}\right]-\frac{1}{m}\nabla_{X}|X|^{2}\\
        &=\int_{M^n}|\nabla^2 h|^2+|\nabla Y|^2-2Ric(\nabla h, Y)+\frac{2}{m}div(X^{\flat}\otimes X^{\flat})(X)\\
        &\qquad +n(n-2)(\sigma_{k}-\lambda)^{2}-(n-2)\lambda\frac{1}{m}|X|^{2}+\frac{(n-2)}{m}\sigma_{k}|X|^{2}-\frac{1}{m}\nabla_{X}|X|^{2}.\\
    \end{split}
\end{equation}
Again, the same argument based on Lemma 2.1 of \cite{petersen2009rigidity}, allow us to deduce that
\begin{equation}\label{dddd3}
\int_{M^n}Ric(\nabla h,\nabla h)dv_{g}=\int_{M^n}n^{2}(\sigma_{k}-\lambda)^2+\frac{2n}{m}(\sigma_{k}-\lambda)|X|^{2}+\frac{1}{m^{2}}|X|^{4}-|\nabla^2 h|^2.
\end{equation}
Replacing \eqref{ddd1}, \eqref{dddd4} and \eqref{dddd3} back into \eqref{tt1}, we get

\begin{equation}\label{ddeedd}
    \begin{split}\int_{M^n}2Ric(\nabla h,Y)dv_{g}&=\int_{M^n}2|\nabla^2 h|^2+\frac{2}{m}div(X^{\flat}\otimes X^{\flat})(X-Y)-2n(\sigma_{k}-\lambda)^{2}\\
        &\qquad -\frac{(n+2)}{m}(\sigma_{k}-\lambda)|X|^{2}-\frac{1}{m}\nabla_{X}|X|^{2}-\frac{1}{m^{2}}|X|^{4}.
    \end{split}
\end{equation}
Since 
\begin{equation*}
    \begin{split}div((X^{\flat}\otimes X^{\flat})(\nabla h))&=div(X^{\flat}\otimes X^{\flat})(\nabla h)+g(\nabla^{2}h,X^{\flat}\otimes X^{\flat})\\
    &=div X.g(X,\nabla h)+g(\nabla_{X}X,\nabla h)+g(\nabla^{2}h,X^{\flat}\otimes X^{\flat}),
    \end{split}
\end{equation*}
we have that \eqref{ddeedd} can be rewritten as follows

\begin{equation}\label{ddeeddee}
    \begin{split}\int_{M^n}2Ric(\nabla h,Y)dv_{g}&=\int_{M^n}2|\nabla^2 h|^2-\frac{2}{m}g(\nabla^{2}h,X^{\flat}\otimes X^{\flat})-2n(\sigma_{k}-\lambda)^{2}\\&\qquad-\frac{(n+2)}{m}(\sigma_{k}-\lambda)|X|^{2}-\frac{1}{m}\nabla_{X}|X|^{2}-\frac{1}{m^{2}}|X|^{4}.
    \end{split}
\end{equation}
Combining \eqref{norma} with \eqref{ddeeddee} we produce the desired result.

\end{myproof}

\begin{myproof}{Theorem}{\ref{T8}}
Since the Hodge-de Rham decomposition is
orthogonal on $L^{2}(M)$, we get
\begin{equation*}
    \int_{M^n}g(\nabla h, X)dv_{g}=\int_{M^n}g(\nabla h, \nabla h+Y)dv_{g}=\int_{M^n}|\nabla h|^{2}dv_{g}.
\end{equation*}
Therefore, if
\begin{equation*}
    \int_{M^n}g(\nabla h, X)dv_{g}\leq0,
\end{equation*}
we obtain that $\nabla h=0$  and, consequently, $X=Y$. Now, since $Y$ is a divergence free vector field, we deduce
\[0=div Y=div X=n(\sigma_{k}-\lambda)+\frac{1}{m}|X|^{2},\]
which implies that 
\[n(\sigma_{k}-\lambda)=-\frac{1}{m}|X|^{2}.\]
On the other hand, from the fundamental equation \eqref{fundamental2*}, we have
\begin{equation*}
        \begin{split}\langle\nabla_{X}X,X\rangle&=(\sigma_{k}-\lambda)|X|^{2}+\frac{1}{m}|X|^{4}=-\frac{1}{nm}|X|^{4}+\frac{1}{m}|X|^{4}=\frac{n-1}{nm}|X|^{4}.
        \end{split}
        \end{equation*}
Note that  
 \begin{equation*}
        \begin{split}div (|X|^{2}X)&=|X|^{2} div X+\langle\nabla |X|^{2},X\rangle\\
        &=|X|^{2} div X+\nabla_{X}|X|^{2}\\
        &=|X|^{2} div X+2\langle \nabla_{X}X,X\rangle\\
         &=2\langle \nabla_{X}X,X\rangle.
        \end{split}
        \end{equation*}
Hence
\[0=\int_{M}div(|X|^{2}X)=2\int_{M}\frac{n-1}{nm}|X|^{4}.\]
Therefore $X=0$ and $\sigma_{k}=\lambda$.

\end{myproof}

\begin{myproof}{Theorem}{\ref{T4}}
    As we already know, the fundamental equation
    \[
       \nabla^{2}f-\frac{1}{m}df\otimes df=(\sigma_{k}-\lambda)g,
    \]
    leads to
    \begin{equation}\label{trace}
        L(f):=e^{\frac{f}{m}}\divergence(e^{-\frac{f}{m}}\nabla f)=\Delta f-\frac{1}{m}|\nabla f|^2=n(\sigma_{k}-\lambda),
    \end{equation}
    and because we suppose that $ L(f)\geqslant0 $ we must then admit
    that $ \sigma_{k}-\lambda\geqslant0 $. So, if we now
    take a cut-off function $\psi:M\rightarrow\mathbb{R}$ satisfying
    \begin{equation*}
        0\leqslant\psi\leqslant1\mbox{ on } M,\quad
        \psi\equiv 1\hspace{0,2cm}\text{in}\hspace{0,2cm} B(r),\quad
        \textup{supp}{(\psi)}\subset B(2r)\quad \text{and}\quad
        |\nabla\psi|\leqslant\frac{K}{r},
    \end{equation*}
    where $ K>0 $ is a real constant, we are in place to conclude that
    \begin{equation*}
        \begin{split}
            n\int_{B(r)}
            \left(
            \sigma_{k}-\lambda
            \right)
            d\mu
            &=
            \int_{B(r)}
            n\psi \left(
            \sigma_{k}-\lambda
            \right)
            d\mu\leqslant
            \int_{B(2r)}
            n\psi \left(
            \sigma_{k}-\lambda
            \right)
            d\mu=
            \int_{B(2r)}\psi L(f)d\mu=
            \\
            &=
            -\int_{B(2r)}g(\nabla\psi,\nabla f)d\mu
            \leqslant
            \int_{B(2r))\setminus B(r)}|\nabla\psi||\nabla f|d\mu
            \leqslant  \\
            &
             \leqslant K\int_{B(2r)\setminus B(r)}\frac{|\nabla f|}{r}d\mu
            \leqslant
            2K\int_{M\setminus B(r)}\dfrac{|\nabla f|}{d(x,x_0)}d\mu,
        \end{split}
    \end{equation*}
    from what it follows that
    \begin{align*}
        0\leqslant
        \int_M
        \left(
        \sigma_{k}-\lambda
        \right)
        d\mu
        &=
        \lim_{r\to\infty}
        \int_{B(r)}
        \left(
        \sigma_{k}-\lambda
        \right)
        d\mu\leqslant
        \frac{2K}{n}\lim_{r\to\infty}
        \int_{M\setminus B(r)}\frac{|\nabla f|}{d(x,x_0)}d\mu=0.
    \end{align*}
    Henceforth, we have that $
   L(f)=\sigma_{k}-\lambda=0 $ which proves the
    theorem.
    
\end{myproof}

\section{Invariant by translation examples of quasi $k$-Yamabe solitons}\label{sectioninvariant}

In  this  section, we provide a method to construct  examples  of conformally flat quasi $k$-Yamabe gradient solitons.  We focus our attention on solitons whose solutions are invariant under the action of translation group. More precisely, we consider the Riemannian metric
\begin{equation*}
\delta=\sum_{i=1}^{n}dx_{i}\otimes dx_{i},
\end{equation*}
in coordinates $x=(x_{1},\dots,x_{n})$ of $\mathbb{R}^{n}$, where $n\geq 3$. For an arbitrary choice of non zero vector $\alpha=(\alpha_{1},\dots,\alpha_{n})$ we define the translation function $\xi:\mathbb{R}^{n}\rightarrow\mathbb{R}$  by
\begin{equation*}
\xi(x_{1},\dots,x_{n})=\alpha_{1}x_{1}+\dots+\alpha_{n}x_{n}.
\end{equation*}

Next, we consider that $\mathbb{R}^{n}$ admits a group of symmetries consisting of translations \cite{olver2000applications} and we then look for smooth functions $\varphi, f:(a,b)\subset\mathbb{R}\rightarrow\mathbb{R}$, with $\varphi>0$, such that the compositions 
\[f=f\circ\xi:\xi^{-1}(a,b)\longrightarrow\mathbb{R},\qquad \varphi=\varphi\circ\xi: \xi^{-1}(a,b)\longrightarrow\mathbb{R},\]
satisfies 
\begin{equation}\label{invarianciafundamental}\nabla^{2}f-\frac{1}{m}df\otimes df=(\sigma_{k}-\lambda)g\quad \text{with}\quad  g=\varphi^{-2}\delta
\end{equation}


What has been said above is summed up in the next result.

\begin{theorem}\label{invarianciath}Let $(\mathbb{R}^{n},\delta)$ and $f=f\circ\xi$, $ \varphi=\varphi\circ\xi$ as above. Then $\delta_{ij}\varphi^{-2}$ is a quasi $k$-Yamabe gradient soliton metric with potential function $f=-m\log u$ if, and only if, 

	\begin{equation}\label{eq:01}
	u''+2\frac{u'\varphi'}{\varphi}=0,
	\end{equation}
	\begin{equation}\label{eq:02}
	b_{n,k}\left[k\varphi\varphi''-\frac{n}{2}(\varphi')^{2}\right](\varphi')^{2(k-1)}||\alpha||^{2k}-m\varphi\varphi'\frac{u'}{u}||\alpha||^{2}=\lambda,
	\end{equation}
	where 
	\begin{equation*}b_{n,k}=\frac{(n-1)!}{k!(n-k)!}(-1)^{k-1}\frac{1}{2^{k-1}}.
	\end{equation*}

\end{theorem}

\begin{myproof}{Theorem}{\ref{invarianciath}}
	It is well known that for the conformal metric  $\bar{g}=\varphi^{-2}\delta$, the Ricci curvature is given by \cite{besse2007einstein}:	
\begin{equation}\label{Ricci}Ric_{\bar{g}}=\frac{1}{\varphi^{2}}\Big{\{}(n-2)\varphi \nabla^{2}_{\delta}\varphi+[\varphi\Delta_{\delta}\varphi-(n-1)|\nabla_{\delta}\varphi|^{2}]\delta\Big{\}},
\end{equation}
and consequently, the scalar curvature $R_{\bar{g}}$ on conformal metric
	is given by
	\begin{equation}\label{escalar conforme}R_{\bar{g}}=(n-1)(2\varphi\Delta_{\delta}\varphi-n|\nabla_{\delta}\varphi|^{2}).
	\end{equation}
	
In order to compute the Schouten Tensor on the conformal geometry $A_{\bar{g}}$ we evoke the expression
	\begin{equation*}
A_{\bar{g}}=\frac{1}{n-2}\left[Ric_{\bar{g}}-\frac{R_{\bar{g}}}{2(n-1)}\bar{g}\right].
\end{equation*}
Therefore, from \eqref{Ricci} and \eqref{escalar conforme} we deduce that
		\begin{equation*}A_{\bar{g}}=\frac{\nabla^{2}_{\delta}\varphi}{\varphi}-\frac{|\nabla_{\delta}\varphi|^{2}}{2\varphi^{2}}\delta.
	\end{equation*}

Denote by $\varphi_{x_{i}}$, $\varphi_{x_{i},x_{j}}$ the derivative of $\varphi$ with respect the variables $x_i$ and $x_{i}x_{j}$, respectively. That being said, since we are assuming that $\varphi(\xi)$ and $f(\xi)$ are functions of $\xi=\alpha_{1}x_{1}+\dots+\alpha_{n}x_{n}$, we get

	\begin{equation*}
	\varphi_{,x_{i}}=\varphi'\alpha_{i},\qquad  f_{,x_{i}}=f'\alpha_{i},\qquad \hspace{0.2cm}\varphi_{,x_{i}x_{j}}=\varphi''\alpha_{i}\alpha_{j},\qquad f_{,x_{i}x_{j}}=f''\alpha_{i}\alpha_{j}.\hspace{0.2cm}
	\end{equation*}
Hence

	\begin{equation*}(\bar{g}^{-1}A_{\bar{g}})_{ij}=\varphi\varphi''\alpha_{i}\alpha_{j}-\frac{1}{2}(\varphi')^{2}||\alpha||^{2}\delta_{ij}.
	\end{equation*}
The eigenvalues of $\bar{g}^{-1}A_{\bar{g}}$ are $\theta=-\frac{1}{2}(\varphi')^{2}||\alpha||^{2}$ with multiplicity $(n-1)$, and $\mu=(\varphi\varphi''-\frac{1}{2}(\varphi')^{2})||\alpha||^{2}$ with multiplicity $1$. The formula for $\sigma_{k}$ can be found easily by the binomial expansion of $(x-\theta)^{n-1}(x-\mu)$

		\begin{eqnarray}\label{k curvature} \sigma_{k}&=&\frac{(n-1)!}{k!(n-k)!}\left[(n-k)\theta+k\mu\right]\theta^{k-1}\nonumber\\
	&=&\frac{(n-1)!}{k!(n-k)!}(-1)^{k-1}\frac{1}{2^{k-1}}\left[k\varphi\varphi''-\frac{n}{2}(\varphi')^{2}\right](\varphi')^{2(k-1)}||\alpha||^{2k}.
	\end{eqnarray}

	Now, in order to compute the Hessian of $u$ relatively to $\bar{g}$ we evoke the expression
	\begin{equation*}
	(\nabla_{\bar{g}}^{2}u)_{ij}=u_{x_{i},x_{j}}-\sum_{k=1}^{n}\Gamma_{ij}^{k}u_{x_{k}},
	\end{equation*}
	where the Christoffel symbol $\Gamma_{ij}^{k}$ for distinct $i,j,k$ are given by
	\begin{equation}\Gamma_{ij}^{k}=0,\ \Gamma_{ij}^{i}=-\frac{\varphi_{x_{j}}}{\varphi},\ \Gamma_{ii}^{k}=\frac{\varphi_{x_{k}}}{\varphi}\;\ \mbox{and}\;\ \Gamma_{ii}^{i}=-\frac{\varphi_{x_{i}}}{\varphi}.\nonumber
	\end{equation}
	Therefore,
		\begin{eqnarray}\label{hessian}(\nabla_{\bar{g}}^{2}u)_{ij}&=&u_{x_{i},x_{j}}+\varphi^{-1}(\varphi_{x_{i}}u_{,x_{j}}+\varphi_{x_{j}}u_{x_{i}})-\delta_{ij}\sum_{k}\varphi^{-1}\varphi_{x_{k}}u_{x_{k}}\nonumber\\
	&=&\alpha_{i}\alpha_{j}u''+(2\alpha_{i}\alpha_{j}-\delta_{ij}||\alpha||^{2})\varphi^{-1}\varphi'u'.
	\end{eqnarray}

If we make the change $u=e^{-\frac{f}{m}}$, then the fundamental soliton equation in \eqref{invarianciafundamental} can be rewritten as follows:
    \begin{equation}\label{abcd}\nabla^2 u=-\frac{u}{m}(\sigma_{k}-\lambda)g.
    \end{equation}
Substituting \eqref{k curvature} and \eqref{hessian} into \eqref{abcd} and considering $i\neq j$ we obtain

	\begin{equation*}
	\alpha_{i}\alpha_{j}\left(u''+2\frac{\varphi'u'}{\varphi}\right)=0.
	\end{equation*}
If there exist $i,j$, $i\neq j$ such that $\alpha_{i}\alpha_{j}\neq
	0$, then we get
	\begin{equation*}u''+2\frac{u'\varphi'}{\varphi}=0,
	\end{equation*}
which provides equation \eqref{eq:01}. For $i=j$, substituting \eqref{k curvature} and \eqref{hessian} into \eqref{abcd} we obtain \eqref{eq:02}. 

Now, we need to consider the case $\alpha_{k_{0}}=1$, $\alpha_{k}=0$ for $k\neq k_{0}$. In this case, substituting \eqref{hessian} into \eqref{abcd} we obtain
\begin{equation}\label{opop1}
-\frac{u}{m}(\sigma_k-\lambda)\frac{1}{\varphi^2}=-\frac{\varphi'u'}{\varphi},
\end{equation}
for $i\neq k_{0}$, that is, $\alpha_{i}=0$ with $i=j$, and
\begin{equation}\label{opop2}
-\frac{u}{m}(\sigma_k-\lambda)\frac{1}{\varphi^2}=u''+(2-1)\frac{\varphi'u'}{\varphi},
\end{equation}
for $i=k_{0}$, that is, $\alpha_{k_{0}}=1$ with $i=j=k_{0}$.

However, \eqref{opop1} and \eqref{opop2} are equivalent to equations \eqref{eq:01} and \eqref{eq:02}. This completes the demonstration.

\end{myproof}

In what follows, we provide examples illustrating Theorem \ref{invarianciath}. We pointed out that Example \ref{examplee1}, Example \ref{Example noncompact} and  Example \ref{examplee2} are contemplated in the next examples.

\begin{example}
In Theorem \ref{invarianciath} consider an arbitrary direction  $\alpha=(\alpha_{1},\alpha_{2},\dots,\alpha_{n})$. Therefore the functions
\[f(x_{1},\dots,x_{n})=-m\log\left[k_{1}(\alpha_{1}x_{1}+\dots+\alpha_{n}x_{n})\right],\qquad \varphi(x_{1},\dots, x_{n})=k_{0},\qquad k_{0}, k_{1}\in(0,\infty),\]
provide a family of quasi $k$-Yamabe gradient soliton on $\{(x_{1},\dots, x_{n})\in\mathbb{R}^n\; ; \;  \alpha_{1}x_{1}+\dots+\alpha_{n}x_{n}>0\}$ with soliton constant $\lambda=0$.
\end{example}

\begin{example}
In Theorem \ref{invarianciath} consider an arbitrary direction $\alpha=(\alpha_{1},\alpha_{2},\dots,\alpha_{n})$. Therefore the functions
\[f(x_{1},\dots,x_{n})=-m\log\left(\frac{k_{1}}{\alpha_{1}x_{1}+\dots+\alpha_{n}x_{n}}\right),\quad \varphi(x_{1},\dots, x_{n})=k_{0}(\alpha_{1}x_{1}+\dots+\alpha_{n}x_{n}),\]
where $k_{0}, k_{1}\in(0,\infty)$ provide a family of quasi $k$-Yamabe gradient soliton on the Euclidean subset $\{(x_{1},\dots, x_{n})\in\mathbb{R}^n\; ;\;  \alpha_{1}x_{1}+\dots+\alpha_{n}x_{n}>0\}$ with soliton constant 
\[\lambda=-\frac{n!}{k!(n-k)!}(-1)^{k-1}\left(\frac{k_{0}^{2}||\alpha||^{2}}{2}\right)^{k}+mk_{0}^{2}||\alpha||^{2}.\]
\end{example}

\begin{example}\label{101010}
In Theorem \ref{invarianciath} consider the direction $\alpha=(0,0,\dots, 0, 1)$. Therefore the functions
\[f(x_{1},\dots,x_{n})=-m\log\left(\frac{k_{1}}{x_{n}}\right),\qquad \varphi(x_{1},\dots, x_{n})=k_{0}x_{n},\quad k_{0}, k_{1}\in(0,\infty),\]
provide a family of quasi $k$-Yamabe gradient soliton on $\{(x_{1},\dots, x_{n})\in\mathbb{R}^n\; ;\;  x_{n}>0\}$ with soliton constant \[\lambda=-\frac{n!}{k!(n-k)!}(-1)^{k-1}\left(\frac{k_{0}^{2}}{2}\right)^{k}+mk_{0}^{2}.\]
\end{example}


\begin{example}
In Theorem \ref{invarianciath} consider an arbitrary direction $\alpha=(\alpha_{1},\alpha_{2},\dots,\alpha_{2k})$. Therefore the functions
\[f(x_{1},\dots,x_{2k})=k_{0}\qquad \varphi(x_{1},\dots, x_{2k})=k_{1}e^{\alpha_{1}x_{1}+\cdots+\alpha_{2k}x_{2k}},\quad k_{0}, k_{1}\in(0,\infty),\]
provide a family of quasi $k$-Yamabe gradient soliton on $\mathbb{R}^{2k}$ with soliton constant $\lambda=0$.
\end{example}

\begin{example}
In Theorem \ref{invarianciath} consider an arbitrary direction $\alpha=(\alpha_{1},\alpha_{2},\dots,\alpha_{n})$. Therefore the functions
\[f(x_{1},\dots,x_{n})=k_{0}\qquad \varphi(x_{1},\dots, x_{n})=k_{2}\{k_{1}k+(2k-n)\left[\alpha_{1}x_{1}+\cdots+\alpha_{n}x_{n}\right]\}^{\frac{2k}{2k-n}},\quad k_{0}, k_{1}, k_{2}\in(0,\infty),\]
provide a family of quasi $k$-Yamabe gradient soliton on $\mathbb{R}^{n}$ with soliton constant $\lambda=0$.
\end{example}

\end{document}